# GEOMETRIC QUANTIZATION BY PATHS

## PART II: CONSTRUCTION OF THE PREQUANTUM GROUPOIDS

PATRICK IGLESIAS-ZEMMOUR

ABSTRACT. In Part I, we constructed the prequantum groupoid for connected and simply connected parasymplectic spaces. This second part extends the construction to *arbitrary connected parasymplectic diffeological spaces* $(X, \omega)$. The presence of a non-trivial fundamental group $\pi_1(X)$ introduces an obstruction to the additivity of the naive integration function on the space of loops, measured by a surfacic 2-cocycle.

We resolve this obstruction by disarticulating the periods into a canonical three-tiered geometric hierarchy: spherical periods, toric periods, and surfacic periods. The *Total Group of Periods* $P_\omega \subset \mathbf{R}$ absorbs the parasymplectic areas of the algebraic relations in $\pi_1(X)$. A choice of a rectifying cochain $c$ for the surfacic coboundary yields a family of prequantum groupoids $\mathbf{T}_c$, each obtained as a quotient of the space of paths by a coherence relation $\sim_c$. We prove that each $\mathbf{T}_c$ is a fibrating diffeological groupoid, that its isotropy group is strictly connected and isomorphic to the total torus of periods $T_\omega = \mathbf{R}/P_\omega$, and that it carries a unique left- and right-invariant prequantum 1-form $\boldsymbol{\lambda}_c$ whose curvature is $\mathrm{trg}^*\omega - \mathrm{src}^*\omega$. The family $\{\mathbf{T}_c\}_c$ is an affine space directed by the character group $\mathrm{Hom}(\pi_1(X), T_\omega)$.

The analysis of the symmetry group of these groupoids, and their comparison with the classical prequantum bundles of Kostant–Souriau, is deferred to a subsequent work.

## INTRODUCTION

This paper is the continuation of [PIZ25a], in the framework of diffeology [PIZ13]. In Part I, we constructed the prequantum groupoid $\mathbf{T}_\omega$ for connected and simply connected parasymplectic spaces $(X, \omega)$. The construction proceeds by distilling the space of paths Paths(X), equipped with the Action Form $\mathsf{K}\omega$, into a diffeological groupoid whose space of morphisms is the quotient of Paths(X) by the

*Date*: October 2, 2026.

2020 *Mathematics Subject Classification.* Primary 53D50, 58A05; Secondary 22A22, 55R65, 58A10.

*Key words and phrases.* Diffeology, Geometric Quantization, Prequantum Groupoid, Closed 2-forms, Parasymplectic Spaces, Paths, Loops.

The author thanks the Hebrew University of Jerusalem, Israel, for his continuous academic support. The TEX source of this paper is available at the GitHub PIZ-Diffeology Archive: `https://github.com/p-i-z/Diffeology-Archives`.

relation identifying paths with the same endpoints and the same action modulo periods.

The assumption of simple connectedness is restrictive. Many physically relevant systems, such as the torus or surfaces of higher genus, are not simply connected. In these cases, the geometry of the space X introduces new periods —surfacic periods— that are not captured by the spherical or toric cycles alone. A naive application of the Part I construction would fail to produce a connected isotropy group, and the quotient would not be a fibrating groupoid.

In this Part II, we extend the construction to arbitrary connected diffeological spaces. Our strategy relies on a direct path approach: the obstruction to the additivity of the prequantum action is identified directly on the space of loops, and measured by a cocycle $\tau$ on the fundamental group $\pi_1(X, x_0)$. By incorporating the values of $\tau$ on the group relations into the Total Group of Periods $P_\omega$, we ensure that the prequantum structure closes on the quotient torus $T_\omega = \mathbf{R}/P_\omega$. A choice of rectifying cochain $c$ for the surfacic coboundary then produces a prequantum groupoid $\mathbf{T}_c$, and varying $c$ produces a family indexed by the character group $\mathrm{Hom}(\pi_1(X), T_\omega)$.

The remainder of the paper is organized as follows. Section I constructs the three-tiered period hierarchy and the Total Group of Periods. Section II constructs the rectified integration function and the Chasles function $\Phi_c$, whose vanishing defines the coherence relation. Section III defines the prequantum groupoids $\mathbf{T}_c$ and establishes their structural properties. Section IV applies the construction to the surface of genus 2 and to the torus $T^2$. We close with a short conclusion.

We make no assumptions regarding the status of the objects as manifolds, infinite-dimensional manifolds, or Lie groups. All constructions are internal to the category of diffeological spaces, and the operations of taking quotients and forming functional spaces require no topological or regularity conditions to yield well-defined smooth structures.

## I. The Three-Tiered Period Hierarchy

The torus of periods, quotient of the real line $\mathbf{R}$ by the group of periods $P_\omega$, plays a fundamental role in geometric quantization as the universal carrier of the quantum phase. When X is simply connected, the periods originate entirely from smooth 2-spheres. In the general non-simply connected case, the periods undergo a canonical three-tiered filtration: spherical periods, toric periods, and surfacic periods.

We construct the total torus of periods through a sequence of successive diffeological subductions:

$$\mathbf{R} \twoheadrightarrow T_{\mathrm{sph}} = \mathbf{R}/P_{\mathrm{sph}} \twoheadrightarrow T_{\mathrm{tor}} = \mathbf{R}/P_{\mathrm{tor}} \twoheadrightarrow T_\omega = \mathbf{R}/P_\omega .$$

**1. The Spherical Periods.**

Let $(X, \omega)$ be a connected parasymplectic diffeological space. Recall the chain-homotopy operator $\mathsf{K}\colon \Omega^k(X) \to \Omega^{k-1}(\mathrm{Paths}(X))$ and the Action Form $\mathsf{K}\omega \in \Omega^1(\mathrm{Paths}(X))$.

We denote by $L_0 = \mathrm{Loops}_0(X, x_0)$ the connected component of contractible loops based at a chosen basepoint $x_0 \in X$. The restriction of $\mathsf{K}\omega$ to $L_0$ is closed. The *Group of Spherical Periods* $P_{\mathrm{sph}}$ is the subgroup of $\mathbf{R}$ defined by integrating $\mathsf{K}\omega$ over all closed cycles in $L_0$:

$$P_{\mathrm{sph}} = \left\{ \int_\sigma \mathsf{K}\omega \;\middle|\; \sigma \in \mathrm{Loops}(L_0) \right\} \subset \mathbf{R}.$$

Geometrically, each cycle $\sigma \in \mathrm{Loops}(L_0)$ represents a smooth 2-sphere in X, and $\int_\sigma \mathsf{K}\omega = \int_{S^2} \omega$. We assume that $P_{\mathrm{sph}}$ is a strict subgroup ($P_{\mathrm{sph}} \neq \mathbf{R}$). The quotient $T_{\mathrm{sph}} = \mathbf{R}/P_{\mathrm{sph}}$ is the *Torus of Spherical Periods*. It carries a canonical closed 1-form $\theta_{\mathrm{sph}} \in \Omega^1(T_{\mathrm{sph}})$ characterized by $\pi_{\mathrm{sph}}^*(\theta_{\mathrm{sph}}) = dt$, where $\pi_{\mathrm{sph}}\colon \mathbf{R} \to T_{\mathrm{sph}}$ is the canonical projection.

**2. The Toric Periods.**

When $\pi_1(X) \neq 0$, the space of free loops decomposes into its connected components:

$$\mathrm{Loops}(X) = \coprod_{i \in \pi_0(\mathrm{Loops}(X))} L_i.$$

On each component $L_i$, the restriction of the closed 1-form $\mathsf{K}\omega$ defines a local period group:

$$P_i = \left\{ \int_\sigma \mathsf{K}\omega \;\middle|\; \sigma \in \mathrm{Loops}(L_i) \right\} \subset \mathbf{R}.$$

Geometrically, a cycle $\sigma \in \mathrm{Loops}(L_i)$ corresponds to a smooth map from a 2-torus into X, $\Sigma\colon T^2 \to X$, and $\int_\sigma \mathsf{K}\omega = \int_{T^2} \Sigma^*(\omega)$.

The *Group of Toric Periods* $P_{\mathrm{tor}}$ is the subgroup of $\mathbf{R}$ generated by the union of all local periods:

$$P_{\mathrm{tor}} = \left\langle \bigcup_{i \in \pi_0(\mathrm{Loops}(X))} P_i \right\rangle.$$

Since the contractible component $L_0$ is one of these components, we have the canonical inclusion $P_{\mathrm{sph}} \subset P_{\mathrm{tor}} \subset \mathbf{R}$. We assume that $P_{\mathrm{tor}}$ is a strict subgroup. The quotient $T_{\mathrm{tor}} = \mathbf{R}/P_{\mathrm{tor}}$ is the *Torus of Toric Periods*. It is equipped with its canonical closed 1-form $\theta_{\mathrm{tor}}$ characterized by $\pi_{\mathrm{tor}}^*(\theta_{\mathrm{tor}}) = dt$.

**3. The Group of Indices.**

To structure the homotopy of X, we fix a basepoint $x_0 \in X$ and consider the space of based loops $\mathrm{Loops}(X, x_0)$. The set of connected components $\mathscr{I} = \pi_0(\mathrm{Loops}(X, x_0))$ is canonically identified with the fundamental group $\pi_1(X, x_0)$. Concatenation of based loops induces a group multiplication on $\mathscr{I}$.

We choose a *basis of loops* $\mathscr{B} = \{\ell_i\}_{i \in \mathscr{I}}$ such that:

(1) For each index $i \in \mathscr{I}$, the loop $\ell_i \in \mathrm{Loops}(X, x_0)$ belongs to the connected component $L_i$.
(2) For the identity element $0 \in \mathscr{I}$, we choose the constant loop $\ell_0 = \boldsymbol{x_0}$.

**4. The Intermediate Integration Function.**

We define the intermediate integration function $\phi\colon \mathrm{Loops}(X, x_0) \to T_{\mathrm{tor}}$ component-wise relative to the chosen basis $\mathscr{B}$. For any loop $\ell$ in the component $L_i$ (indexed by $i \in \mathscr{I}$):

$$\phi(\ell) = \pi_{\mathrm{tor}}\left(\int_{\ell_i}^{\ell} \mathsf{K}\omega\right) \in T_{\mathrm{tor}},$$

where $\pi_{\mathrm{tor}}\colon \mathbf{R} \to T_{\mathrm{tor}}$ is the canonical projection, and the integral is evaluated along any smooth path in the component $L_i$ connecting the basis loop $\ell_i$ to $\ell$. Because any two paths connecting $\ell_i$ to $\ell$ within $L_i$ differ by a closed cycle in $\mathrm{Loops}(L_i)$, their integrals differ by an element of the local period group $P_i \subset P_{\mathrm{tor}}$. Therefore, $\phi(\ell)$ is strictly well-defined in $T_{\mathrm{tor}}$, independent of the path chosen within $L_i$. By construction, $\phi$ satisfies the differential identity $\phi^*(\theta_{\mathrm{tor}}) = \mathsf{K}\omega \restriction \mathrm{Loops}(X, x_0)$.

**5. The Surfacic Cocycle.**

The function $\phi$ measures the parasymplectic action relative to the chosen basis $\mathscr{B}$. However, $\phi$ is generally not a homomorphism with respect to loop concatenation. Consider two based loops $\ell \in L_i$ and $\ell' \in L_j$. Their concatenation $\ell \vee \ell'$ belongs to the component $L_{i\cdot j}$. The integral defining $\phi(\ell \vee \ell')$ decomposes geometrically into two steps: first from the basis loop $\ell_{i\cdot j}$ to the concatenation $\ell_i \vee \ell_j$, and then from $\ell_i \vee \ell_j$ to $\ell \vee \ell'$. Using the additivity of path integration over concatenated homotopies established in Part I, the second term splits into the individual integrals modulo $P_{\mathrm{tor}}$. Consequently, we obtain the fundamental relation:

$$\phi(\ell \vee \ell') = \phi(\ell) + \phi(\ell') + \tau(i, j),$$

where the defect of additivity is measured by the *surfacic cocycle* $\tau\colon \mathscr{I} \times \mathscr{I} \to T_{\mathrm{tor}}$:

$$\tau(i, j) = \pi_{\mathrm{tor}}\left(\int_{\ell_{i\cdot j}}^{\ell_i \vee \ell_j} \mathsf{K}\omega\right) \in T_{\mathrm{tor}}.$$

**Proposition (The Surfacic Cocycle).** The map $\tau\colon \mathscr{I} \times \mathscr{I} \to T_{\mathrm{tor}}$ is a group 2-cocycle of $\mathscr{I}$ with values in the abelian group $T_{\mathrm{tor}}$. That is, for all $i, j, k \in \mathscr{I}$:

$$\tau(i, j) + \tau(i \cdot j, k) = \tau(j, k) + \tau(i, j \cdot k).$$

**Note (Geometric Meaning).** The cocycle $\tau$ admits a direct interpretation in terms of parasymplectic areas. If two elements $i, j \in \mathscr{I}$ commute, the difference $\tau(i, j) - \tau(j, i)$ is the area of the 2-torus swept out by sliding $\ell_i$ past $\ell_j$. If $i$ and $j$ do not commute, the area is captured by the accumulated cocycle over a full relation in the fundamental group.

*Proof.* Let $i, j, k \in \mathscr{I}$. We evaluate the action of deforming $\ell_{i\cdot j\cdot k}$ to the triple concatenation $\ell_i \vee \ell_j \vee \ell_k$ along two distinct homotopic decomposition paths in $\mathrm{Loops}(\mathrm{X}, x_0)$.

*1. Left-parenthesized decomposition.* Decompose the deformation from $\ell_{(i\cdot j)\cdot k}$ to $(\ell_i \vee \ell_j) \vee \ell_k$ into two steps: first $\ell_{(i\cdot j)\cdot k} \to \ell_{i\cdot j} \vee \ell_k$, contributing $\tau(i\cdot j, k)$ modulo $\mathrm{P_{tor}}$; then $\ell_{i\cdot j} \vee \ell_k \to (\ell_i \vee \ell_j) \vee \ell_k$, contributing $\tau(i, j)$ by additivity. The total is $\tau(i, j) + \tau(i\cdot j, k)$.

*2. Right-parenthesized decomposition.* Similarly, the deformation from $\ell_{i\cdot(j\cdot k)}$ to $\ell_i \vee (\ell_j \vee \ell_k)$ through $\ell_i \vee \ell_{j\cdot k}$ contributes $\tau(i, j\cdot k) + \tau(j, k)$.

*3. Associativity homotopy.* The two endpoints $(\ell_i \vee \ell_j) \vee \ell_k$ and $\ell_i \vee (\ell_j \vee \ell_k)$ are connected by the canonical smooth reparametrization $\alpha\colon [0,1] \to \mathrm{Loops}(\mathrm{X}, x_0)$. The associated 2-plot $\Phi(s, t) = \alpha(s)(t)$ factors through a 1-plot, so its pullback of $\omega$ vanishes identically. Hence $\int_\alpha \mathbf{K}\omega = 0$.

*4. Conclusion.* The closed cycle formed by the first path, the associativity homotopy, and the reverse of the second path lies in $\mathrm{L}_{i\cdot j\cdot k}$, so its total integral is zero modulo $\mathrm{P_{tor}}$. Hence $\tau(i, j) + \tau(i\cdot j, k) = \tau(j, k) + \tau(i, j\cdot k)$ in $\mathrm{T_{tor}}$. □

**6. The Central Extension.**

The 2-cocycle $\tau$ defines a group extension $\tilde{\Gamma}$ of the fundamental group $\mathscr{I}$ by the torus of toric periods $\mathrm{T_{tor}}$. As a set, $\tilde{\Gamma} = \mathscr{I} \times \mathrm{T_{tor}}$, equipped with the multiplication law:

$$(i, u)\cdot(j, v) = (i\cdot j,\ u + v + \tau(i, j)).$$

This extension is strictly central: the fiber elements $(0, a)$ commute with all elements of $\tilde{\Gamma}$.

**7. The Total Group of Periods.**

Let $\mathrm{F}(\mathscr{I})$ be the free group generated by $\mathscr{I}$, and let $\mathscr{R} \subset \mathrm{F}(\mathscr{I})$ be the normal subgroup of relations, so that $\mathscr{I} \simeq \mathrm{F}(\mathscr{I})/\mathscr{R}$. Let $\sigma\colon \mathscr{I} \to \tilde{\Gamma}$ be the canonical section $\sigma(i) = (i, 0)$. By the universal property of free groups, $\sigma$ extends uniquely to a group homomorphism $\Psi\colon \mathrm{F}(\mathscr{I}) \to \tilde{\Gamma}$. Because $\Psi$ lifts the canonical projection, it maps $\mathscr{R}$ into the central fiber $\mathrm{T_{tor}}$.

For $w \in \mathscr{R}$, the *accumulated cocycle* $\mathscr{T}(w) \in \mathrm{T_{tor}}$ is defined by $\Psi(w) = (0, \mathscr{T}(w))$.

**Definition Total Group of Periods.** Let $\mathrm{H_{surf}} \subset \mathrm{T_{tor}}$ be the subgroup generated by the accumulated cocycle values on all relations:

$$\mathrm{H_{surf}} = \left\langle \left\{ \mathscr{T}(w) \mid w \in \mathscr{R} \right\} \right\rangle.$$

The *Total Group of Periods* is the preimage of $\mathrm{H_{surf}}$ in $\mathbf{R}$ under $\pi_{\mathrm{tor}}$:

$$\mathrm{P}_\omega = \pi_{\mathrm{tor}}^{-1}(\mathrm{H_{surf}}) \subset \mathbf{R}.$$

We assume $P_\omega \neq \mathbf{R}$. The quotient $T_\omega = \mathbf{R}/P_\omega \simeq T_{\rm tor}/H_{\rm surf}$ is the *Total Torus of Periods.*

**Note (Algebraic Significance).** By construction, the projection of $\tau$ to $T_\omega$ is a coboundary: $[\pi_\omega \circ \tau] = 0 \in H^2(\pi_1(X), T_\omega)$. Equivalently, there exists a 1-cochain $c\colon \mathscr{I} \to T_\omega$ such that:

$$\pi_\omega(\tau(i,j)) = c(i) + c(j) - c(i \cdot j).$$

This cochain is the datum that rectifies the integration function, as explained in Section II.

## II. Rectified Integration and Chasles Function

With the Total Torus of Periods $T_\omega$ established, we construct the rectified integration function on the based loop space. Because the surfacic cocycle $\tau$ becomes a coboundary when projected to $T_\omega$, the naive integration can be rectified into a genuine homomorphism with respect to loop concatenation. The rectified function is then extended to free loops, and used to define the Chasles function on pairs of paths with common endpoints.

### 8. The Rectified Integration on Based Loops.

By the construction of the Total Group of Periods $P_\omega$, the projection of the surfacic cocycle to $T_\omega$ is a coboundary. Therefore, there exists a 1-cochain $c\colon \mathscr{I} \to T_\omega$ such that:

$$\pi_\omega(\tau(i,j)) = c(i) + c(j) - c(i \cdot j),$$

for all $i, j \in \mathscr{I}$. We fix such a cochain once and for all, normalized by $c(0) = 0$ for the identity element.

Recall the integration function $\phi\colon \mathrm{Loops}(X, x_0) \to T_{\rm tor}$, satisfying the identity $\phi(\ell \vee \ell') = \phi(\ell) + \phi(\ell') + \tau(i,j)$ for $\ell \in L_i$ and $\ell' \in L_j$. Project this function on the Total Torus of Periods and let $\phi' = \pi_\omega \circ \phi$. Then:

$$\begin{aligned}\phi'(\ell \vee \ell') &= \phi'(\ell) + \phi'(\ell') + \pi_\omega(\tau(i,j)) \\ &= \phi'(\ell) + \phi'(\ell') + c(i) + c(j) - c(i \cdot j).\end{aligned}$$

That is:

$$\phi'(\ell \vee \ell') + c(i \cdot j) = \big(\phi'(\ell) + c(i)\big) + \big(\phi'(\ell') + c(j)\big).$$

Defining then

$$\phi_c(\ell) = \phi'(\ell) + c(i) \quad \text{for all} \quad \ell \in L_i,$$

we obtain the following proposition.

**Proposition (Properties of $\phi_c$).** The map $\phi_c$ satisfies:

1. **Normalization:** $\phi_c(\boldsymbol{x_0}) = 0$.
2. **Additivity:** for all $\ell, \ell' \in \mathrm{Loops}(X, x_0)$, $\phi_c(\ell \vee \ell') = \phi_c(\ell) + \phi_c(\ell')$.

(3) **Primitive property:** $\phi_c^*(\theta_\omega) = \mathbf{K}\omega \restriction \mathrm{Loops}(X, x_0)$, where $\theta_\omega$ is the canonical 1-form on $T_\omega$.

(4) **Family:** If $c'$ is another cochain satisfying the same coboundary relation, then $\chi = c' - c$ is a character $\mathscr{I} \to T_\omega$, and $\phi_{c'} = \phi_c + \chi \circ \mathrm{comp}$, where $\mathrm{comp}\colon \mathrm{Loops}(X, x_0) \to \mathscr{I}$ is the component map. The family $\{\phi_c\}_c$ is an affine space directed by $\mathrm{Hom}(\mathscr{I}, T_\omega)$.

*Proof.* *(1) Normalization.* The constant loop $\boldsymbol{x_0}$ belongs to $L_0$, and $c(0) = 0$, so $\phi_c(\boldsymbol{x_0}) = 0$.

*(2) Additivity.* Immediate from the identity established above.

*(3) Primitive property.* On each component $L_i$, the cochain $c(i)$ is constant, so its differential vanishes. Hence $\phi_c^*(\theta_\omega) \restriction L_i = \phi'^*(\theta_\omega) \restriction L_i = \phi^*(\theta_{\mathrm{tor}}) \restriction L_i = \mathbf{K}\omega \restriction L_i$. Summing over components gives the result.

*(4) Family.* If $c'$ satisfies the same coboundary relation, then $\chi = c' - c$ satisfies $\chi(i) + \chi(j) - \chi(i \cdot j) = 0$, so $\chi$ is a homomorphism. The relation $\phi_{c'} = \phi_c + \chi \circ \mathrm{comp}$ is immediate from the definition. Freeness and transitivity of the action follow. □

**9. Extension to the Free Loop Space.**

The rectified integration $\phi_c$ is defined on the based loop space. We extend it to the free loop space Loops(X).

Let $\ell \in \mathrm{Loops}(X)$, and let $x = \ell(0) = \ell(1)$ be its basepoint. Choose a smooth path $\delta \in \mathrm{Paths}(X)$ from $x_0$ to $x$. The concatenation $\delta \vee \ell \vee \bar{\delta}$ is a loop based at $x_0$. We define the value of the free loop integration function by:

$$\phi_c(\ell) = \phi_c(\delta \vee \ell \vee \bar{\delta}) \in T_\omega.$$

**Proposition (Independence of the connecting path).** The extension is well-defined and independent of the choice of $\delta$. For this reason we denote it by the same symbol $\phi_c$.

*Proof.* Let $\delta'$ be another path from $x_0$ to $x$. The two loops $\delta' \vee \ell \vee \bar{\delta}'$ and $C \vee (\delta \vee \ell \vee \bar{\delta}) \vee \bar{C}$, where $C = \delta' \vee \bar{\delta}$, differ by the insertion of a backtrack. The integral of $\mathbf{K}\omega$ along a backtrack vanishes, because it amounts to integrating the 2-form $\omega$ on a 1-plot. Since $\phi_c$ is additive and $T_\omega$ is abelian, the two values coincide. The extension is unique: it is the only map $\mathrm{Loops}(X) \to T_\omega$ whose restriction to $\mathrm{Loops}(X, x_0)$ is $\phi_c$ and which is additive under concatenation. □

**10. The Chasles Function.**

We recall the space of pairs of paths with coincident endpoints:

$$\mathrm{Paths}_2(X) = \{(\gamma, \gamma') \in \mathrm{Paths}(X) \times \mathrm{Paths}(X) \mid \mathrm{ends}(\gamma) = \mathrm{ends}(\gamma')\}.$$

Let $(\gamma, \gamma') \in \mathrm{Paths}_2(\mathrm{X})$, so that $\gamma(0) = \gamma'(0) = x$ and $\gamma(1) = \gamma'(1) = x'$. The concatenation $\gamma \vee \bar{\gamma}'$ forms a closed loop based at $x$.

We define the *Chasles function* $\Phi_c\colon\ \mathrm{Paths}_2(\mathrm{X}) \to \mathrm{T}_\omega$ by:

$$\Phi_c(\gamma, \gamma') = \phi_c(\gamma \vee \bar{\gamma}'),$$

where $\phi_c$ is the free loop integration function of Article 9. The definition involves no auxiliary path, because $\phi_c$ is defined directly on the free loop space.

**Proposition (Chasles Identity).** For any three paths $\gamma, \gamma', \gamma''$ with the same endpoints:

$$\Phi_c(\gamma, \gamma') + \Phi_c(\gamma', \gamma'') = \Phi_c(\gamma, \gamma'').$$

**Remark (Family).** If $c'$ is another cochain and $\chi = c' - c$, then $\Phi_{c'}(\gamma, \gamma') = \Phi_c(\gamma, \gamma') + \chi(\gamma \vee \bar{\gamma}')$. The family $\{\Phi_c\}_c$ is an affine space directed by $\mathrm{Hom}(\mathscr{I}, \mathrm{T}_\omega)$.

*Proof.* The loop $\gamma \vee \bar{\gamma}''$ and the concatenation $(\gamma \vee \bar{\gamma}') \vee (\gamma' \vee \bar{\gamma}'')$ differ by a backtrack, whose integral vanishes. By additivity of $\phi_c$ on free loops sharing a basepoint, $\Phi_c(\gamma, \gamma'') = \Phi_c(\gamma, \gamma') + \Phi_c(\gamma', \gamma'')$. ☐

## III. The Prequantum Groupoids

All the geometric ingredients are now in place to construct the prequantum groupoids for an arbitrary connected parasymplectic space $(\mathrm{X}, \omega)$. We distill the Primordial Quantum Broth $(\mathrm{Paths}(\mathrm{X}), \mathsf{K}\omega)$ into a family of prequantum groupoids $\{\mathbf{T}_c\}_c$, one for each choice of cochain $c$ resolving the surfacic coboundary.

### 11. The Quantum Broth.

Before crystallizing the structure into a groupoid, we recognize the fundamental object from which all quantum geometry derives. We call the pair $(\mathrm{Paths}(\mathrm{X}), \mathsf{K}\omega)$ the *Primordial Quantum Broth.*[1] This infinite-dimensional space contains every possible transition of the system. By integration, the Action Form $\mathsf{K}\omega$ assigns a potential phase to every finite fluctuation of a motion. However, this space is noisy: it contains infinitely many distinct paths that are physically indistinguishable because they differ by variations carrying no parasymplectic action modulo periods.

The construction of the prequantum groupoids is the process of *distilling* this broth. We identify paths that carry the same phase information, collapsing the infinite redundancy of the broth into the precise geometry of the groupoid.

[1]This structure has a concrete existence in diffeology; it is not merely a figure of speech. Diffeology was expressly developed to handle such infinite-dimensional structures —a functional space equipped with a genuine differential form— with simplicity and rigor.

**12. The Equivalence Relation.**

For each cochain $c$ resolving the surfacic coboundary, we define the *coherence relation* $\sim_c$ on Paths(X):

$$\gamma \sim_c \gamma' \quad \text{if and only if} \quad \mathrm{ends}(\gamma) = \mathrm{ends}(\gamma') \quad \text{and} \quad \Phi_c(\gamma, \gamma') = 0 \in T_\omega.$$

Here, $\mathrm{ends}(\gamma) = (\gamma(0), \gamma(1)) \in X \times X$ is the ends map, and $\Phi_c\colon \mathrm{Paths}_2(X) \to T_\omega$ is the Chasles function of Article 10.

This relation identifies paths sharing the same endpoints whose relative parasymplectic action, corrected by the cochain $c$ on the components of the loop space, vanishes in the Total Torus of Periods.

**13. Definition of the Groupoids.**

Given a resolution $c$ of the surfacic cocycle $\tau$, the quotient of the space of paths by the coherence relation $\sim_c$ defines the morphisms of the prequantum groupoid.

**Definition (The Prequantum Groupoids).** The *Prequantum Groupoid* $\mathbf{T}_c$ is the category internal to diffeology defined by:

$$\mathrm{Obj}(\mathbf{T}_c) = X, \qquad \mathrm{Mor}(\mathbf{T}_c) = \mathcal{Y}_c = \mathrm{Paths}(X)/\sim_c.$$

The source and target maps $\mathrm{src}, \mathrm{trg}\colon \mathcal{Y}_c \to X$ are induced by the end-points evaluation maps $\hat{0}, \hat{1}\colon \mathrm{Paths}(X) \to X$:

$$\mathrm{src}([\gamma]_c) = \gamma(0), \quad \mathrm{trg}([\gamma]_c) = \gamma(1).$$

The composition of composable morphisms is induced by the concatenation of paths:

$$[\gamma]_c \cdot [\gamma']_c = [\gamma \vee \gamma']_c, \quad \text{for} \quad \gamma(1) = \gamma'(0).$$

The unit map $x \mapsto \mathbf{1}_x = [\boldsymbol{x}]_c$ assigns to each classical motion $x \in X$ the class of the constant path at $x$.

**14. Structural Properties.**

The resulting object is a fibrating diffeological groupoid equipped with a prequantum 1-form.

**Theorem (The Prequantum Structure).** Let $(X, \omega)$ be a connected parasymplectic diffeological space with $\omega \neq 0$, and let $c$ be a resolution of the surfacic coboundary. The groupoid $\mathbf{T}_c$ satisfies:

(1) It is a fibrating diffeological groupoid over X: the projection $\mathrm{ends}\colon \mathcal{Y}_c \to X \times X$ is a subduction.
(2) The isotropy group at any point $x \in X$, $\mathbf{T}_{c,x}$, is strictly connected and isomorphic to the Total Torus of Periods $T_\omega$ via the descent of $\phi_c$.
(3) There exists a unique left- and right-invariant 1-form $\boldsymbol{\lambda}_c \in \Omega^1(\mathcal{Y}_c)$, the *prequantum 1-form*, such that $\mathrm{class}_c^*(\boldsymbol{\lambda}_c) = \mathbf{K}\omega$, where $\mathrm{class}_c\colon \mathrm{Paths}(X) \to \mathcal{Y}_c$ is the quotient projection.

(4) The curvature of $\boldsymbol{\lambda}_c$ satisfies $d\boldsymbol{\lambda}_c = \mathrm{trg}^*\,\omega - \mathrm{src}^*\,\omega$.

*Proof. 1. Groupoid structure and subduction.* The additivity of $\phi_c$ guarantees that $\sim_c$ is compatible with concatenation. The ends map factors through the quotient, and since $\mathrm{ends}\colon \mathrm{Paths}(X) \to X \times X$ is a subduction, its factorization on $\mathcal{Y}_c$ is also a subduction.

*2. Isotropy.* The isotropy at $x_0$ consists of the morphism classes of loops in $\mathrm{Loops}(X, x_0)$. Two classes coincide iff $\Phi_c(\ell, \ell') = 0$, which by additivity of $\phi_c$ is equivalent to $\phi_c(\ell) = \phi_c(\ell')$. Hence the induced map $\overline{\phi_c}\colon \mathbf{T}_{c,x_0} \to \mathrm{im}(\phi_c)$ is well-defined and injective. That $\mathrm{im}(\phi_c) = T_\omega$ follows from the same argument as in the simply connected case, using the contractible component $L_0$ and the fact that the path-primitive $\mathbf{F}$ on $L_0$ has image all of $\mathbf{R}$. The extension to arbitrary $x$ uses conjugation by a path from $x$ to $x_0$, which preserves $\sim_c$ because $\Phi_c$ is conjugation-invariant.

*3. The prequantum form.* Apply the diffeological descent criterion to $\mathsf{K}\omega$ and the subduction $\mathrm{class}_c$. If $\mathrm{class}_c \circ P = \mathrm{class}_c \circ P'$ for two plots $P, P'$, then for all $r$, $\mathrm{ends}(P(r)) = \mathrm{ends}(P'(r))$ and $\Phi_c(P(r), P'(r)) = 0$. Hence $\int_{P'(r)}^{P(r)} \mathsf{K}\omega \in P_\omega$, which is locally constant because $P_\omega$ is discrete. The exterior derivative vanishes, so $\mathsf{K}\omega(P) = \mathsf{K}\omega(P')$ on each connected domain, and the descent criterion applies.

*4. Curvature.* By the chain-homotopy identity, $d\mathsf{K}\omega = \hat{1}^*\omega - \hat{0}^*\omega$. Pulling back along $\mathrm{class}_c$ and using the relation $\hat{1}_{\mathrm{Paths}} = \mathrm{trg} \circ \mathrm{class}_c$, $\hat{0}_{\mathrm{Paths}} = \mathrm{src} \circ \mathrm{class}_c$, we get $\mathrm{class}_c^*(d\boldsymbol{\lambda}_c) = \mathrm{class}_c^*(\mathrm{trg}^*\,\omega - \mathrm{src}^*\,\omega)$. Since $\mathrm{class}_c$ is a subduction, $d\boldsymbol{\lambda}_c = \mathrm{trg}^*\,\omega - \mathrm{src}^*\,\omega$. □

**Remark (The family).** The construction yields one prequantum groupoid for each cochain $c$ resolving the surfacic coboundary. The family $\{\mathbf{T}_c\}_c$ is parametrized by the character group $\mathrm{Hom}(\pi_1(X), T_\omega)$ once a base cochain is fixed, and the whole family is an affine space directed by that group. The study of the symmetry group of these groupoids, the comparison with the classical prequantum bundles of Kostant and Souriau, and the analysis of the moduli space of the family, are the subject of a subsequent work.

## IV. Examples

We illustrate the construction with two archetypal examples: a surface of genus 2, where the quantum phase is generated purely by the surfacic periods, and the torus $T^2$, where the three tiers collapse into one.

### 15. Explicit Construction: The Surface of Genus 2.

Let $\Sigma_2$ be a compact, orientable surface of genus 2, equipped with an area form $\omega$. This example illustrates the mechanics of the Total Group of Periods, because the spherical and toric periods vanish completely, leaving only the surfacic periods to generate the quantum phase.

*1. Vanishing of spherical and toric periods.* The surface $\Sigma_2$ has the hyperbolic plane $\mathbf{H}^2$ as its contractible universal cover, so $\pi_2(\Sigma_2) = 0$ and $P_{sph} = \{0\}$. Every smooth map $T^2 \to \Sigma_2$ has degree zero (because every abelian subgroup of a surface group of genus at least two is cyclic), so $P_{tor} = \{0\}$ and $T_{tor} = \mathbf{R}$.

*2. Surfacic period.* The fundamental group $\pi_1(\Sigma_2, x_0)$ is generated by four loops $a_1, b_1, a_2, b_2$ with the single defining relation $R = [a_1, b_1] \cdot [a_2, b_2] = 0$. Geometrically, the relation corresponds to the concatenated loop $\gamma_R$, whose lift to $\mathbf{H}^2$ is the boundary of the fundamental octagon. The accumulated cocycle is:

$$\mathscr{T}(R) = \int_{x_0}^{\gamma_R} \mathsf{K}\omega = \iint_{\text{Octagon}} \tilde{\omega} = \text{Area}(\Sigma_2) = A.$$

Hence $P_\omega = A\mathbf{Z}$ and $T_\omega = \mathbf{R}/A\mathbf{Z}$.

*3. Contractible loops sweep out* $\mathbf{R}$. A loop $\ell \in L_0$ lifts to a closed loop $\tilde{\ell}$ in $\mathbf{H}^2$ bounding a disk D. Since $\mathbf{H}^2$ has infinite area, there is a smooth 1-parameter family of loops whose bounded areas cover all of $\mathbf{R}$. Projecting back to $\Sigma_2$ and to $T_\omega$, the image of $L_0$ under $\phi_c$ covers the whole torus.

*4. Restoring the fundamental group relation.* The concatenation of the eight generator loops yields $\gamma_R$, whose class must equal the identity morphism $\mathbf{1}_{x_0}$ for the groupoid to faithfully represent $\pi_1(\Sigma_2)$. The phase difference between $\gamma_R$ and the constant loop is A. Quotienting by $P_\omega = A\mathbf{Z}$ forces $\phi_c(\gamma_R) \equiv 0$, identifying the octagon boundary with the constant loop.

**16. The Torus** $T^2$**.**

Let $X = \mathbf{R}^2/\mathbf{Z}^2$ be the standard torus equipped with the symplectic area form $\omega = dx \wedge dy$. The fundamental group is $\pi_1(X) \simeq \mathbf{Z}^2$, generated by the two canonical loops A and B. The toric periods are $P_{tor} = \mathbf{Z}$.

The fundamental relation of the torus is the commutator $ABA^{-1}B^{-1} = 1$, and the surfacic period associated with it is the area of the fundamental domain, which is 1. Since $1 \in \mathbf{Z} = P_{tor}$, the surfacic contribution does not enlarge the period group. Therefore:

$$P_\omega = P_{tor} = \mathbf{Z}.$$

The prequantum isotropy is $\mathbf{T}_{c,x} \simeq \mathbf{R}/\mathbf{Z} \simeq S^1$, recovering the standard prequantization of the torus. This contrasts with higher genus surfaces, where the surfacic period is strictly larger than the toric period.

## Conclusion

This paper has constructed the prequantum groupoids $\mathbf{T}_c$ for arbitrary connected parasymplectic diffeological spaces. The construction proceeds by disarticulating the periods of the parasymplectic form into a three-tiered hierarchy (spherical, toric, and surfacic), incorporating the surfacic cocycle into the Total Group of Periods, rectifying the naive integration function by a choice of cochain $c$, and quotienting the space of paths by the resulting coherence relation. Each resulting groupoid is fibrating, its isotropy group is the total torus of periods, and it carries a canonical prequantum 1-form whose curvature is the difference of the parasymplectic form on the target and source.

The family of prequantum groupoids is indexed by the character group of the fundamental group, and the whole family forms an affine space directed by this group. The analysis of the symmetry group of these groupoids, the identification with the classical prequantum bundles of Kostant and Souriau, and the study of the moduli space of the family, will be presented in subsequent work.

Einstein Institute of Mathematics, The Hebrew University of Jerusalem, Campus Givat Ram, 9190401 Israel.

*Email address*: piz@math.huji.ac.il

*URL*: http://math.huji.ac.il/~piz

*URL*: https://github.com/p-i-z/Diffeology-Archives